\documentstyle[amssymb,amsfonts]{amsart}

\newcommand\dist{{\operatorname{dist}}}
\newcommand\R{{{\mathbf R}}}
\newcommand\Z{{{\mathbf Z}}}
\newenvironment{proof}{\noindent {\bf Proof} }{\endprf\par}
\def \endprf{\hfill  {\vrule height6pt width6pt depth0pt}\medskip}
\def\emph#1{{\it #1}}
\def\textbf#1{{\bf #1}}

\theoremstyle{plain}
  \newtheorem{theorem}[subsection]{Theorem}
  
  \newtheorem{proposition}[subsection]{Proposition}
  \newtheorem{lemma}[subsection]{Lemma}
  \newtheorem{corollary}[subsection]{Corollary}

\theoremstyle{remark}
  \newtheorem{remark}[subsection]{Remark}

\theoremstyle{definition}

\include{psfig}

\begin{document}

\title[Energy-critical NLW in $\R^{1+3}$]{Spacetime bounds for the energy-critical nonlinear wave equation in three spatial dimensions}
\author{Terence Tao}
\address{Department of Mathematics, UCLA, Los Angeles CA 90095-1555}
\email{ tao@@math.ucla.edu}
\subjclass{35L15}
\thanks{The author is supported by a grant from the Packard Foundation.  The author also thanks Kenji Nakanishi, Fabrice Planchon, and Tristan Roy for corrections, references, and other helpful comments.}

\vspace{-0.3in}
\begin{abstract}
Results of Struwe, Grillakis, Struwe-Shatah, Kapitanski, Bahouri-Shatah, Bahouri-G\'erard, and Nakanishi have established global wellposedness, regularity, and scattering in the energy class for the energy-critical nonlinear wave equation $\Box u = u^5$ in $\R^{1+3}$, together with a spacetime bound
$$ \| u \|_{L^4_t L^{12}_x(\R^{1+3})} \leq M(E(u))$$
for some finite quantity $M(E(u))$ depending only on the energy $E(u)$ of $u$.  We reprove this result, and show that this quantity obeys a bound of at most exponential type in the energy, and specifically $M(E) \leq C (1+E)^{C E^{105/2}}$ for some absolute constant $C > 0$.
The argument combines the quantitative local potential energy decay estimates of these previous papers with arguments used by Bourgain and the author for the analogous nonlinear Schr\"odinger equation.
\end{abstract}

\maketitle

\section{Introduction}

In this paper we revisit the energy-critical (quintic) defocusing nonlinear wave equation
\begin{equation}\label{nlw}
\Box u = u^5
\end{equation}
in three spatial dimensions, where $u: I \times \R^3 \to \R$ is a real scalar field (with time restricted to some time interval $I$),
and $\Box u := - \partial_{tt} u + \Delta u$ is the d'Lambertian.  To avoid technicalities let us restrict attention to \emph{classical solutions}
to this equation, by which we mean solutions which are infinitely smooth and are compactly supported in space for each fixed time $t$.  (It turns out that more general classes of finite energy solution can be written as strong limits in the energy class of such classical solutions.)
Then one can easily verify that solutions to these equations have a conserved energy
\begin{equation}\label{energy}
E(u) := \int_{\R^3} \frac{1}{2} |\partial_t u |^2 + \frac{1}{2} |\nabla u|^2 + \frac{1}{6} |u|^6\ dx
\end{equation}
and an invariant scaling
\begin{equation}\label{scaling}
u(t,x) \mapsto \frac{1}{\lambda^{1/2}} u(\frac{t}{\lambda}, \frac{x}{\lambda}).
\end{equation}
Note that the energy is \emph{critical} (scale-invariant).  Informally speaking, this means that the nonlinear components of \eqref{nlw} are significant at all time scales, especially when the energy is large, and so the proper analysis of this equation has proven to be somewhat delicate.
Nevertheless, after an intensive effort by many researchers (\cite{segal-weak}, \cite{rauch}, \cite{pecher}, \cite{struwe}, \cite{g_waveI}, \cite{struweshatah.energy}, \cite{grillakis.semilinear}, \cite{kap.energy}, \cite{gsv}, \cite{bs}, \cite{bg}, \cite{nakscatter}, \cite{nak2}; see \cite{shatah-struwe} or \cite{sogge:wave} for a detailed treatment of the topic) this equation was demonstrated to be well-behaved in the energy class.  In particular, given any classical data there was a unique classical global solution, and given any finite energy data there was a unique global finite energy solution which obeyed\footnote{It remains an open question whether there are any other solutions in the energy class, which are have an infinite $L^4_t L^{12}_x$ norm; we will not address this ``unconditional uniqueness'' issue here, but see \cite{mas} for some recent progress on this question.} the additional spacetime integrability property\footnote{In the literature the $L^5_t L^{10}_x$ norm is often used instead of the $L^4_t L^{12}_x$ norm, but bounds on one easily imply bounds on the other by interpolation and Strichartz estimates.  Also, standard arguments allow one to convert these sorts of spacetime bounds into quantitative persistence of regularity bounds, controlling the $C^0_t H^s_x \cap C^1_t H^{s-1}_x$ norm of the solution by the $H^s_x \times H^{s-1}_x$ norm of the initial data for any $s > 1$.} of lying in the space $L^4_t L^{12}_x(\R \times \R^3)$ (see the next section for definitions of these norms).
It was even demonstrated in \cite{bg}, \cite{nak2} that this global $L^4_t L^{12}_x(\R \times \R^3)$ norm was bounded uniformly by some unspecified quantity depending on the energy.  The purpose of this note is to explicitly establish this bound.

\begin{theorem}\label{main}  Let $(u_0,u_1) \in \dot H^1_x(\R^3) \times L^2_x(\R^3)$ be initial data with the energy bound
$$ \int_{\R^3} \frac{1}{2} |u_1 |^2 + \frac{1}{2} |\nabla u_0|^2 + \frac{1}{6} |u_0|^6\ dx \leq E.$$
Then there exists a unique global solution $u \in C^0_t \dot H^1_x(\R \times \R^3) \cap \dot C^1_t L^2_x(\R \times \R^3) \cap L^4_t L^{12}_x(\R \times \R^3)$ with the spacetime bound
\begin{equation}\label{412-main} \| u \|_{L^4_t L^{12}_x(\R \times \R^3)} \leq C (1 + E)^{C E^{105/2}}
\end{equation}
for some absolute constant $C > 0$.
\end{theorem}

\begin{remark} The only novelty of this theorem is the bound \eqref{412-main} in the case of large energy; the remaining portions of the theorem were already obtained in previous works \cite{struwe}, \cite{g_waveI}, \cite{struweshatah.energy}, \cite{grillakis.semilinear}, \cite{kap.energy}, \cite{bs}, \cite{bg}, \cite{nakscatter}, \cite{nak2}, whereas the small energy theory follows immediately from Strichartz estimates (see \eqref{strichartz}), with results in that case dating back to \cite{rauch}, \cite{pecher}. In particular, in \cite{bg} it was shown via a concentration compactness method that the $L^4_t L^{12}_x(\R \times \R^3)$ norm of $u$ was bounded by some unspecified function $M(E)$ of the energy, though as remarked in \cite{bg} it was not clear exactly what \emph{explicit} bound $M(E)$ could be extracted
from that argument; it seems unlikely that the bound will end up being of exponential type of better, and is more likely to be tower-exponential in nature.  The arguments in \cite{nak2} (which also extend to the Klein-Gordon equation) also in principle yield an explicit bound on $M(E)$.  This argument uses the induction on energy method of Bourgain, and would thus would be expected to yield tower-exponential bounds, though in the massless case $m=0$ an exponential bound can be extracted (Nakanishi, personal communication). 
One can also obtain other scale-invariant bounds on the solution $u$ by further application of Strichartz estimates, since
the nonlinearity $u^5$ will now enjoy bounds in $L^1_t L^2_x(\R \times \R^3)$, but we shall not detail these here.  As observed in \cite{bg}, \cite{bs}, these types of bounds also imply scattering results for \eqref{nlw} in the energy class; see \cite{bg}, \cite{bs} for further details.
Indeed the arguments here could be used to provide a (somewhat complicated) reproof of several of the results cited above.  
\end{remark}

The arguments are broadly in the same spirit as (though not quite identical to) a similar spacetime bound for spherically symmetric solutions of the energy-critical defocusing nonlinear \emph{Schr\"odinger} equation in three dimensions, see \cite{tao}.  The main idea is to isolate those ``bubbles'' in of spacetime where the $L^4_t L^{12}_x$ norm (as well as the mass and energy) is concentrating, and then apply the energy decay arguments used in previous papers to control the $L^4_t L^{12}_x$ norm in the domain of influence of the most concentrated such ``bubble''.  One then removes this bubble, thus decrementing the energy, and then repeats the procedure.  The argument thus falls under the general class of \emph{induction on energy} arguments introduced by Bourgain \cite{bourg.critical}, although the induction here is fairly simple compared to other instances of that argument and in particular will not lead to tower-exponential bounds or worse here.  The methods also share some similarities with those in \cite{nak2}, although the spaces used, the exact Morawetz estimate used, and the combinatorial devices used to organise the various concentrations of the solution in spacetime are all slightly different.

The methods used are fairly standard; besides the induction on energy method, and the energy decay estimate (Proposition \ref{g}) which also plays a key role in earlier work, we rely primarily on Strichartz estimates, energy estimates, 
H\"older's inequality, finite speed of propagation, and some elementary combinatorics (such as the pigeonhole principle, and the exclusion of exceptional intervals as in \cite{tao}) to organise the various regions of spacetime\footnote{In fact, almost all of the combinatorics shall be confined to organising the time variable.}  in which energy concentration has been detected.  There are 
however two slightly exotic tools that are worth mentioning here.  Firstly
there is an \emph{inverse Sobolev theorem} (Lemma \ref{globinvsob}, which shows that whenever the potential energy is large compared to the kinetic energy, then the mass and energy of the solution must be concentrating in a ball.  (A similar result appears in \cite[Lemma 3.5]{bg}, for a related purpose.)  Secondly, one can exploit bounds on the energy flux via the fundamental solution \eqref{fundsoln} in a simple fashion to control the long-range effects of the nonlinearity, thus obtaining a type of asymptotic stability result for linear evolutions of the solutions; see Corollary \ref{ass}.  This means that except for certain ``exceptional'' regions of spacetime where the scattering solutions have a significant presence, the solution is sustained almost entirely by \emph{short-range} interactions and thus cannot concentrate in too isolated a region of spacetime; this will ultimately lead to a key ``non-lacunarity''
estimate (Lemma \ref{nonlacun}) which is crucial to the argument.  This idea is inspired by a somewhat similar one in \cite{tao}, in which the long-range effects of the nonlinearity were shown to enjoy a certain H\"older regularity and thus could not contribute to isolated concentrations
in spacetime.  These arguments also appear to be related to the high-frequency approximation results in \cite{bg}, though the author was not able to
draw a precise connection.

\begin{remark} In the case when the solution is spherically symmetric one can obtain much better results \cite{gsv}.  We sketch the main details here.  Firstly, a standard Morawetz inequality argument (based on inspecting the time derivative of the outgoing momentum $\int_{\R^3} \partial_t u(t,x) \frac{x}{|x|}\cdot \nabla u(t,x)\ dx$) gives the spacetime bound
$$ \int_\R \int_{\R^3} \frac{|u|^6}{|x|}\ dx dt \lesssim E.$$
Secondly, in the spherically symmetric case one has the radial Sobolev inequality
$$ |u(t,x)| |x|^{1/2} \lesssim E^{1/2}.$$
Combining these two we obtain the $L^8_{t,x}$ bound
$$ \| u \|_{L^8_{t,x}(\R \times \R^3)} \lesssim E^{1/4}.$$
By subdividing the time axis $\R$ into $O(E^6)$ intervals where the $L^8$ norm has size $cE^{-1/2}$ for some small $c$, and then using Strichartz estimates \eqref{strichartz} and energy conservation \eqref{energy-bound-kinetic} on each piece, and then summing up, one can obtain a bound of the form $O(E^2)$ for the $L^4_tL^{12}_x$ norm.
\end{remark}

\begin{remark} There are three independent sources in our argument for the exponential growth.  One arises from a pigeonholing argument required to obtain flux decay (Corollary \ref{space-decay}).  Another arises from a pigeonholing argument required to properly eliminate a bubble of energy concentration from a solution (see the arguments after \eqref{massb}).  The third arises from the induction on energy argument used to deal with the portions of the solution away from the bubbles.  It does not seem likely that one could eliminate all three sources of growth from the argument to obtain a more natural bound (such as a polynomial bound), though one may tentatively conjecture that such a bound exists, as is already indicated in the spherically symmetric case discussed above. On the other hand, the exponent of $105/2$ can certainly be lowered somewhat.
\end{remark}

\begin{remark} There should be no difficulty extending this result to higher dimensions, after adjusting all the numerology; results of this nature already appear in \cite{nak2}.  In five and higher dimensions the nonlinearity is not smooth and so one cannot work exclusively with classical solutions (unless one mollifies the nonlinearity slightly), but this is likely to only be a minor technical nuisance.  The proof of the inverse Sobolev theorem given here relies on the fact that the Sobolev exponent $2^* = 6$ is an even integer, but there is no question that an analogue of the theorem exists in more general dimensions (for instance
by modifying the proof of \cite[Lemma 3.5]{bg}).
\end{remark}

\section{Notation}

We use $X \lesssim Y$ or $X \leq O(Y)$ to denote the estimate $X \leq CY$ for some absolute constant $C$ (which can vary from line to line).

We use $L^q_t L^r_x$ to denote the spacetime norm
$$ \| u \|_{L^q_t L^r_x(\R \times \R^3)} := (\int_\R (\int_{\R^3} |u(t,x)|^r\ dx)^{q/r}\ dt)^{1/q},$$
with the usual modifications when $q$ or $r$ is equal to infinity, or when the domain $\R \times \R^3$
is replaced by a smaller region of spacetime such as $I \times \R^3$.  

We shall have occasional need of the following Littlewood-Paley projection operators.
Let $\varphi(\xi)$ be a bump function adapted to the ball $\{ \xi \in
\R^3: |\xi| \leq 2 \}$ which equals 1 on the ball $\{ \xi \in \R^3:
|\xi| \leq 1 \}$.  Define a \emph{dyadic number} to be any number $N \in 2^\Z$ of the form
$N = 2^j$ where $j \in \Z$ is an integer.  For each dyadic number $N$,
we define the Fourier multipliers
\begin{align*}
\widehat{P_{\leq N} f}(\xi) &:= \varphi(\xi/N) \hat f(\xi)\\
\widehat{P_{> N} f}(\xi) &:= (1 - \varphi(\xi/N)) \hat f(\xi)\\
\widehat{P_N f}(\xi) &:= (\varphi(\xi/N) - \varphi(2\xi/N)) \hat f(\xi).
\end{align*}
We record a specific instance of the \emph{Bernstein inequality}, namely
\begin{equation}\label{bernstein}
\| P_{\leq N} f \|_{L^{12}_x(\R^3)} \lesssim N^{1/4} \| f \|_{L^6_x(\R^3)}.
\end{equation}
This follows readily from Young's inequality and an inspection of the convolution kernel of $P_{\leq N}$.

\section{Fundamental solution, Strichartz, and Sobolev}

In this section we recall some basic estimates concerning the free wave equation and the Sobolev inequality.

Let $u: I \times \R^3 \to \R$ be a classical solution to the equation $\Box u = F$ for some time interval $I$ and
$F: I \times \R^3 \to \R$.  For any time $t_0 \in I$,
define the \emph{free development $u_{t_0}$ to $u$ from time $t_0$} to be the unique solution to the free wave equation $\Box u_{t_0} = 0$ which agrees with $u$ at time $t=t_0$ to first order, thus $u_{t_0}(t_0) = u(t_0)$ and $\partial_t u_{t_0}(t_0) = \partial_t u(t_0)$.  We recall the \emph{Duhamel formula}
\begin{equation}\label{duhamel}
u(t) = u_{t_0}(t) + \int_{t_0}^t \frac{\sin( (t-t') \sqrt{-\Delta} )}{\sqrt{-\Delta}} F(t')\ dt'
\end{equation}
for all $t_0, t \in I$, where we adopt the convention that $\int_{t_0}^t = - \int_t^{t_0}$ when $t < t_0$.  We can write the operator $\frac{\sin( (t-t') \sqrt{-\Delta} )}{\sqrt{-\Delta}}$ more explicitly for $t \neq t'$ as
\begin{equation}\label{fundsoln}
 [\frac{\sin( (t-t') \sqrt{-\Delta} )}{\sqrt{-\Delta}} F(t')](x) = \frac{1}{4\pi(t'-t)} \int_{|x'-x| = |t'-t|} F(t',x')\ dS(x')
 \end{equation}
for all $x \in \R^3$, where $dS$ denotes the surface element.

We also recall the \emph{Strichartz estimate}
\begin{equation}\label{strichartz}
\begin{split}
\| u \|_{L^3_t L^{18}_x(I \times \R^3)} &+ 
\| u \|_{L^4_t L^{12}_x(I \times \R^3)} + 
\| u \|_{L^\infty_t L^6_x(I \times \R^3)} + 
\| \nabla_{t,x} u \|_{L^\infty_t L^2_x(I \times \R^3)}\\ 
&\lesssim \| \nabla_{t,x} u(t_0) \|_{L^2_x(\R^3)} + \| F \|_{L^1_t L^2_x(I \times \R^3)}.
\end{split}
\end{equation}
For a proof, see \cite{ginebre:summarywave}, \cite{sogge:wave}, or \cite{Keel-Tao}.  Many other such estimates are available, but these are the only ones we shall need.

The Sobolev inequality
\begin{equation}\label{sobolev}
\| f \|_{L^6_x(\R^3)} \lesssim \| f \|_{\dot H^1_x(\R^3)}
\end{equation}
is of course very well known.  Now we investigate the circumstances in which this inequality becomes close to being an equality.  (Note that similar considerations also appeared in \cite{bg}.)

\begin{lemma}[Inverse Sobolev inequality]\label{globinvsob} Let $f \in \dot H^1_x(\R^3)$ be such that $\| f \|_{\dot H^1_x(\R^3)} \lesssim E^{1/2}$ and $\|f\|_{L^6_x(\R^3)} \gtrsim \eta$ for some $E, \eta > 0$.  
Then there exists a ball $B(x,r)$ in $\R^3$ such that we have the mass concentration estimate
\begin{equation}\label{br}
 \int_{B(x,r)} |f(y)|^2\ dy \gtrsim E^{-1/2} \eta^{3} r^2.
\end{equation}
and the energy concentration estimate
\begin{equation}\label{br-energy}
 \int_{B(x,r)} |\nabla f(y)|^2\ dy \gtrsim E^{-1/2} \eta^{3}.
 \end{equation}
If in addition we have the stronger estimate $\|P_{\geq N} f \|_{L^6_x(\R^3)} \gtrsim \eta$ for some $N > 0$, then we can take the radius $r$ to be $O(1/N)$.
\end{lemma}

\begin{remark} In the converse direction, if \eqref{br} holds, then an easy application of H\"older's inequality gives
$\|f\|_{L^6_x(\R^3)} \gtrsim E^{-1/4} \eta^{3/2}$.  Thus this theorem is efficient up to polynomial losses.  
\end{remark}

\begin{proof}  By dividing $\eta$ and $f$ by $E^{1/2}$ we may normalise $E=1$; we may also take $f$ to be real-valued.  
We use Littlewood-Paley theory.  Expanding $f = \sum_N P_N f$ and using the triangle inequality, we have
$$\sum_{N_1, \ldots, N_6} |\int_{\R^3} P_{N_1} f \ldots P_{N_6} f| \gtrsim \eta^6.$$
By ordering the $N_j$, we thus have
$$\sum_{N_1 \geq N_2 \geq \ldots \geq N_6} |\int_{\R^3} P_{N_1} f \ldots P_{N_6} f| \gtrsim \eta^6.$$
Observe that we may take $N_2 \sim N_1$ otherwise the integral vanishes.
Using H\"older's inequality, we obtain
$$\sum_{N_1 \sim N_2 \geq \ldots \geq N_6} \| P_{N_1} f \|_{L^2_x(\R^3)} \| P_{N_2} f \|_{L^2_x(\R^3)} \| P_{N_3} f \|_{L^\infty_x(\R^3)} \ldots \| P_{N_6} f \|_{L^\infty_x(\R^3)} \gtrsim \eta^6$$
and hence
$$\sum_{N_1 \sim N_2} \| P_{N_1} f \|_{L^2_x(\R^3)} \| P_{N_2} f \|_{L^2_x(\R^3)} (\sum_{N_3 \leq N_2} \| P_{N_3} f \|_{L^\infty_x(\R^3)})^4 \gtrsim \eta^6.$$
On the other hand, a simple application of Plancherel and Cauchy-Schwarz gives
$$\sum_{N_1 \sim N_2} N_2^2 \| P_{N_1} f \|_{L^2_x(\R^3)} \| P_{N_2} f \|_{L^2_x(\R^3)} \lesssim \| f\|_{\dot H^1_x(\R^3)}^2 \lesssim 1$$
and thus by the pigeonhole principle there exists $N_2$ such that
$$ (\sum_{N_3 \leq N_2} \| P_{N_3} f \|_{L^\infty_x(\R^3)})^4 \gtrsim \eta^6 N_2^{2}.$$
Taking fourth roots and applying the pigeonhole principle again, there exists $N_3$ such that
$$ \| P_{N_3} f \|_{L^\infty_x(\R^3)} \gtrsim \eta^{3/2} N_3^{1/2}.$$
Let us now introduce a bump function $\phi$ supported on $B(0,1)$ whose Fourier transform has magnitude $\sim 1$ on the ball $B(0,100)$,
and write $\psi_{N_3}(x) := N_3^3 (\Delta \phi)(N_3 x)$.  Elementary Fourier analysis shows that
$P_{N_3} f = \tilde P_{N_3} (f * \psi_{N_3})$, where $\tilde P_{N_3}$ is a Fourier multiplier similar to $P_{N_3}$, and in particular is bounded on $L^\infty_x$.  We conclude that there exists an $x$ such that
$$ |f * \psi_{N_3}(x)| \gtrsim \eta^{3/2} N_3^{1/2}$$
or in other words
$$ |\int_{B(x,1/N_3)} f(y) (\Delta \phi)(N_3 (x-y))\ dy| \gtrsim \eta^{3/2} N_3^{-5/2}.$$
The claim \eqref{br} (with $r = 1/N_3)$ now follows from H\"older's inequality, while \eqref{br-energy} follows from one integration by parts
follows by H\"older's inequality.  If we now assume the stronger estimate $\| P_{\geq N} f \|_{L^6_x} \gtrsim \eta$, then the same arguments
lead to
$$ \| P_{N_3} P_{\geq N} f \|_{L^\infty_x(\R^3)} \gtrsim \eta^{3/2} N_3^{1/2}$$
which force $N_3 \gtrsim N$.  One now argues as before (with $P_{N_3}$ replaced by the essentially similar multiplier $P_{N_3} P_{\geq N}$) to obtain
the claim $r = O(1/N)$.
\end{proof}

\section{Energy bounds}

In this section we assume that $u: I \to \R^3$ is a classical solution to \eqref{nlw} with energy $E[u] \leq E$.  As the small energy theory is well understood, we will restrict attention to the case $E \gtrsim 1$.

From the definition of energy we obtain the potential energy bound
\begin{equation}\label{energy-bound-potential}
\| u \|_{L^\infty_t L^6_x(I \times \R^3)} \lesssim E^{1/6}
\end{equation}
and the kinetic energy bound
\begin{equation}\label{energy-bound-kinetic}
\| \nabla_{t,x} u \|_{L^\infty_t L^2_x(I \times \R^3)} \lesssim E^{1/2}.
\end{equation}
From Hardy's inequality we thus conclude that
\begin{equation}\label{hardy}
 \int_{\R^3} \frac{|u(t,y)|^2}{|y-x|^2}\ dy \lesssim E
\end{equation}
for all $t \in I$ and $x \in \R^3$.  
We also have some control on local masses for short times (cf. the corresponding estimates for NLS in
\cite[Section 2.2]{tao}):

\begin{lemma}[Local mass is locally stable]\label{massstable}  Let $B(x,r)$ be a ball, and let $t, t' \in I$.  Then we have
$$ (\int_{B(x,r)} |u(t',y)|^2\ dy)^{1/2} = (\int_{B(x,r)} |u(t,y)|^2\ dy)^{1/2} + O( E^{1/2} |t-t'| ).$$
\end{lemma}

\begin{proof}  We may assume $t' \geq t$. From Minkowski's inequality we have
$$ (\int_{B(x,r)} |u(t',y)|^2\ dy)^{1/2} \leq (\int_{B(x,r)} |u(t,y)|^2\ dy)^{1/2} + \int_t^{t'} (\int_{B(x,r)} |\partial_t u(s,y)|^2\ dy)^{1/2}\ ds
$$
and the claim follows from \eqref{energy-bound-kinetic}.
\end{proof}

\begin{remark} This estimate will combine well with the mass concentration estimate in Lemma \ref{globinvsob}.  Together, they imply that whenever the potential energy is large, there will be a ``bubble'' of spacetime (of spatial and temporal extent $\sim r$) where the solution will be rather large (comparable to $r^{-1/2}$ in some $L^2$ averaged sense).  The task of bounding spacetime norms such as the $L^4_t L^{12}_x$ norm is essentially equivalent to the task of counting how many ``distinct'' bubbles there are in a solution.  In achieving this goal, it will be useful to gather as many upper and lower bounds on the radius and lifespan $r$ of these bubbles as possible.
\end{remark}

Now, let us intersect the spacetime slab $I \times \R^3$ with the forward solid light cone 
$$ \Gamma_+ := \{ (t,x): |x| < t \}$$
defining $\Gamma_+(I) := \Gamma_+ \cap (I \times \R^3)$.  We define the local energies
$$ e(t) := \int_{|x| < t} \frac{1}{2} |\partial_t u(t,x)|^2 + \frac{1}{2} |\nabla u(t,x)|^2 + \frac{1}{6} |u(t,x)|^6\ dx$$
for all $t \geq 0$, thus $0 \leq e(t) \leq E$.  The quantity $e(t)$ is monotone increasing in time, indeed we can easily verify the inequality
$$ \partial_t e(t) \geq \int_{|x| = t} |u(t,x)|^6\ dS(x)$$
where $dS$ denotes the surface element.
Integrating this, we obtain the \emph{energy flux bound}
$$ \int_I \int_{|x| = t} |u(t,x)|^6\ dS(x) dt \lesssim E$$
and then by exploiting spacetime translation invariance and time reversal symmetry we have
\begin{equation}\label{fluxbound}
\int_I \int_{|x'-x| = |t'-t|} |u(t',x')|^6\ dS(x') dt' \lesssim E
\end{equation}
for all $(t',x') \in I \times \R^3$.  We can combine this with the formula \eqref{fundsoln} for the fundamental solution, together with H\"older's inequality, and conclude the useful $L^\infty_x$ bound
\begin{equation}\label{fardecay}
\| \int_{I'} \frac{\sin( (t-t') \sqrt{-\Delta} )}{\sqrt{-\Delta}}(u^5(t')) \|_{L^\infty_x(\R^3)} \lesssim E^{5/6} \dist(t, I')^{-1/2}
\end{equation}
for all subintervals $I' \subseteq I$ and all $t \in I \backslash I'$.  
This bound will allow us to neglect the ``long-range'' effect of the nonlinearity; it is analogous to the similarly subcritical H\"older norm bound in \cite[Lemma 3.4]{tao}, and is used to exclude the possibility of energy concentration in an extremely small ball.
Combining \eqref{fardecay} with Duhamel's formula \eqref{duhamel}, we obtain 

\begin{corollary}[Asymptotic stability]\label{ass}  Let $[t_1,t_2]$ be any time interval in $I$, and let $u_{t_1}, u_{t_2}$ be the free developments of $u$ from times $t_1$ and $t_2$ respectively.  Then for any $t \not \in [t_1,t_2]$ we have
$$ \| u_{t_1}(t) - u_{t_2}(t) \|_{L^\infty_x(\R^3)} \lesssim E^{5/6} \dist(t, [t_1,t_2])^{-1/2}.$$
\end{corollary}

\begin{remark} This result can be compared with the scattering theory in \cite{bs}, which in our language would assert that $u_t$ converges strongly in $C^0_t \dot H^1_x \cap \dot C^1_t L^2_x$ as $t \to +\infty$ or $t \to -\infty$ to an asymptotic state $u_+$ or $u_-$.  Here we establish convergence in a uniform (i.e. $L^\infty_{t,x}$) sense rather than in the energy topology, as long as the time variable is localised to be bounded.
\end{remark}

Now we combine the energy bounds with the Strichartz bounds.  First we observe that largeness of $L^4_t L^{12}_x$ norm implies largeness of potential energy. 

\begin{lemma}[Lower bound on global potential energy]\label{global-pot}  Let $I'$ be a subinterval of $I$ such that
$$ \| u \|_{L^4_t L^{12}_x(I' \times \R^3)} \geq \eta$$
for some $0 < \eta \lesssim E^{1/12}$.  Then there exists $t \in I'$ such that
$$ \| u(t) \|_{L^6_x(\R^3)} \gtrsim E^{-2} \eta^4.$$
\end{lemma}

\begin{proof}  By shrinking $I'$ if necessary we may assume
\begin{equation}\label{412}
 \| u \|_{L^4_t L^{12}_x(I' \times \R^3)} = \eta.
 \end{equation}
From this, the potential energy bound \eqref{energy-bound-potential}, and H\"older, we have
$$ \| u^5 \|_{L^1_t L^2_x( I' \times \R^3 } \lesssim E^{1/6} \eta^4 \lesssim E^{1/2}$$
and hence by Strichartz \eqref{strichartz} and the kinetic energy bound \eqref{energy-bound-kinetic}
\begin{equation}\label{u318}
 \| u \|_{L^3_t L^{18}_x( I' \times \R^3 } \lesssim E^{1/2}.
 \end{equation}
From \eqref{412} and H\"older we conclude that
$$ \| u \|_{L^\infty_t L^6_x( I' \times \R^3) } \gtrsim E^{-2} \eta^4$$
and the claim follows.
\end{proof}

Because of the finite speed of propagation of the fundamental solution \eqref{fundsoln}, we can localise the above argument to the light cone $\Gamma_+$, obtaining

\begin{lemma}[Lower bound on local potential energy]\label{local-pot}  Let $I'$ be a subinterval of $I \cap \R^+$ such that
$$ \| u \|_{L^4_t L^{12}_x(\Gamma_+(I'))} \geq \eta$$
for some $0 < \eta \lesssim E^{1/12}$.  Then there exists $t \in I'$ such that
$$ \int_{|x| \leq t} |u(t,x)|^6 \gtrsim E^{-12} \eta^{24}.$$
\end{lemma}

\begin{proof}  Again we can assume $\| u \|_{L^4_t L^{12}_x(\Gamma_+(I'))} = \eta$.
Let $t_1$ be the upper endpoint of $I'$.  Let $u'$ be the solution to the nonlinear wave equation $\Box u' = (u')^5 1_{\Gamma_+(I')}$ with data $u'(t_1) = u(t_1)$, $\partial_t u'(t_1) = \partial_t u(t_1)$.  Using Strichartz estimates one can show that this solution can be constructed on $I'$, and from finite speed of propagation one can check that the solutions $u$ and $u'$ agree on $\Gamma(I')$.  In particular we see that $\| u' \|_{L^4_t L^{12}_x(\Gamma_+(I'))} = \eta$. Arguing as in the previous lemma we see that $u'$ satisfies \eqref{u318}, and hence
$ \| u \|_{L^3_t L^{18}_x( \Gamma_+(I') } \lesssim E^{1/2}$.  The claim then follows from H\"older's inequality as before.
\end{proof}

We can also combine the argument of Lemma \ref{global-pot} with Lemma \ref{globinvsob} to obtain a mass concentration result also.

\begin{lemma}[Mass concentration estimate]\label{massconc}  Let $I'$ be a subinterval of $I$ such that
$$ \| u \|_{L^4_t L^{12}_x(I' \times \R^3)} \geq \eta$$
for some $0 < \eta \lesssim E^{1/12}$.  Then there exists $t \in I'$ and a ball $B(x,r)$ with $r = O(E^{2/3} \eta^{-4} |I'|)$ such that
we have the mass concentration estimate
$$ \int_{B(x,r)} |u(t,y)|^2\ dy \gtrsim E^{-13/2} \eta^{12} r^2$$
and an energy concentration estimate
$$ \int_{B(x,r)} |\nabla u(t,y)|^2\ dy \gtrsim E^{-13/2} \eta^{12}.$$
\end{lemma}

\begin{proof} Let $N$ be a frequency to be chosen later.  We observe from H\"older's inequality in time, and Bernstein \eqref{bernstein} in space, that
$$ \| P_{\leq N} u \|_{L^4_t L^{12}_x(I' \times \R^3)} \lesssim |I'|^{1/4} N^{1/4} \| u \|_{L^\infty_t L^6_x(I' \times \R^3)}
\lesssim |I'|^{1/4} N^{1/4} E^{1/6}$$
thanks to \eqref{energy-bound-potential}.  Thus by the triangle inequality, we can find $N \sim  E^{-2/3} \eta^4 |I'|^{-1}$ such that
$$ \| P_{>N} u \|_{L^4_t L^{12}_x(I' \times \R^3)} \gtrsim \eta.$$
Applying $P_{>N}$ to \eqref{u318} and then using H\"older we conclude
$$ \| P_{>N} u \|_{L^\infty_t L^6_x( I' \times \R^3) } \gtrsim E^{-2} \eta^4.$$
The claim now follows from Lemma \ref{globinvsob}.
\end{proof}

\begin{remark} We have an upper bound on the size $r$ of the ball on which concentration occurs, but we have no lower bound.  It is possible to establish such a lower bound for ``unexceptional'' intervals $I'$, allowing one to conclude that $r$ is comparable to $|I'|$, up to powers of 
$E$ and $\eta$, by exploiting Strichartz estimates and Corollary \ref{ass}; we will need not do so directly here, because we can work at the scale of the smallest ball of concentration, but see Lemma \ref{nonlacun} or \cite{tao} for some closely related arguments.
\end{remark}

The theory of the energy-critical NLW \eqref{nlw} for large data relies on the following Morawetz-type potential energy decay estimate:

\begin{proposition}[Morawetz-type inequality]\label{g}  Let $[a,b] \subseteq I \cap \R^+$.  Then
$$ \int_{|x| \leq b} |u(b,x)|^6\ dx \lesssim \frac{a}{b} E + (e(b)-e(a)) + (e(b)-e(a))^3.$$
\end{proposition}

Results of this type appear in \cite{g_waveI}, \cite{grillakis.semilinear}, \cite{struwe}, \cite{struweshatah.energy}, \cite{shatah-struwe}, and are based on contracting the stress-energy tensor against various geometrically natural vector fields (such as the scaling vector field $x \cdot \nabla$).  A proof of the above proposition can be found in particular in the appendix of \cite{bg}.  It has the following consequence:

\begin{corollary}[Potential energy decay]  Let $0 < \eta < 1$ and $A \geq E/\eta$.  If we have $[t, A^{CE/\eta} t] \subseteq I \cap \R^+$ for a suitably large absolute constant $C$, then there exists $[t', At'] \subseteq [t, A^{CE/\eta} t]$ such that
$$ \| u \|_{L^\infty_t L^6_x(\Gamma_+([t',At']))} \lesssim \eta.$$
\end{corollary}

\begin{proof} We can locate $\sim CE/\eta$ disjoint intervals of the form $[t'/A,At']$ inside $[t, A^{CE/\eta} t]$.  Since $e(a)$ is an increasing function of $a$, we see from the pigeonhole principle that for one of these intervals, we have $e(b) - e(a) = O(\eta)$ for all
$a,b \in [t'/A,At']$.  The claim then follows by applying Proposition \ref{g} with $a := t'/A$ and $b \in [t',At']$.
\end{proof}

Combining this with with Lemma \ref{local-pot}, we obtain some decay on the $L^4_t L^{12}_x$ norm.

\begin{corollary}[Spacetime norm decay]\label{space-decay}
Let $0 < \eta < 1$ and $A \geq E^3 \eta^{-4}$.  If we have $[t, A^{CE^3 \eta^{-4}} t] \subseteq I \cap \R^+$ for a suitably large absolute constant $C$, then there exists $[t', At'] \subseteq [t, A^{CE^3 \eta^{-4}} t]$ such that
$$ \| u \|_{L^4_t L^{12}_x(\Gamma_+([t',At'])} \lesssim \eta.$$
\end{corollary}

Corollary \ref{space-decay} shows that the cone $\Gamma_+$ contains large slabs on which the $L^4_t L^{12}_x$ norm is negligible.  Unfortunately, it does not directly lead to a uniform bound on the $L^4_t L^{12}_x$ norm on all of $\Gamma_+$, let alone on $I \times \R^3$.  To do this requires some additional arguments, relying on the decay estimate \eqref{fardecay} and the machinery of exceptional and unexceptional intervals from \cite{tao}, as
well as the mass concentration machinery developed above; we turn to this next.

\section{Proof of main theorem}

We are now ready to prove Theorem \ref{main}.
For each energy $E$, let us define $M(E)$ to be the quantity
$$ M(E) := \sup \{ \| u \|_{L^4_t L^{12}_x([t_-,t_+] \times \R^3)} \}$$
where $[t_-,t_+]$ ranges over all compact intervals, and $u$ ranges over all classical solutions to \eqref{nlw} on
the spacetime slab $[t_-,t_+] \times \R^3$ of energy at most $E$.  If $E$ is sufficiently small, one can obtain the bound $M(E) \lesssim 1$
directly from the Strichartz estimate \eqref{strichartz}.  We remark that in \cite{bg} it was shown that $M(E)$ was finite for every $E$, but this was achieved via a compactness method and it is not immediately obvious how to extract a concrete bound for $M(E)$ from that paper; an alternate argument in \cite{nak2}, based on the induction on energy strategy of Bourgain, does in principle give a bound on $M(E)$ which appears to also be exponential in nature, though a precise value for this bound has not been explicitly computed.  Actually, we will not need the results of \cite{bg} or \cite{nak2} here as we shall prove the finiteness of $M(E)$ in the course of our argument.

The local wellposedness theory, described for instance in \cite{shatah-struwe} or \cite{sogge:wave}, allows one to approximate finite energy solutions by classical ones, and to continue them in time, so long as the $L^4_t L^{12}_x$ norm remains bounded.  As such, to 
prove Theorem \ref{main}, it thus suffices by standard limiting arguments  
to establish 
the bound
$$ M(E) \lesssim E^{O( E^{105/2} )}$$
for $E \gtrsim 1$.

We shall use (a rather mild form of) the induction on energy argument of Bourgain \cite{bourg.critical}.  We need some small absolute constants
$$ 1 \gg c_0 \gg c_1 \gg c_2 \gg 0$$
to be chosen later; we always assume that each $c_j$ is chosen to be suitably small depending on the preceding constants.
We also write $C_j := 1/c_j$ for each $j$.  Thus for instance $C_2$ is large compared to $C_0$ and $C_1$.

Since $M$ is already known to be bounded
for small $E$, it will suffice to show

\begin{proposition}[Main inductive step]\label{mainprop}  Let $E \gtrsim 1$.  Then we have the inequality
$$ M(E) \lesssim E^{O(C_2 E^{35})} (1 + M(E - c_0 E^{-33/2})).$$
\end{proposition}

Indeed, one can iterate this proposition about $E^{35/2}$ times and then use the small-energy theory to obtain the desired claim.  Note that it is important here that each $M()$ depends in a linear fashion on the preceding $M()$, otherwise we would obtain tower-exponential bounds as in \cite{bourg.critical} (or even worse, as in \cite{gopher}).

The rest of this section is devoted to the proof of this proposition.  Roughly speaking, the idea is to use the inverse Sobolev theorem (Lemma \ref{globinvsob}, or more precisely Lemma \ref{massconc}) to locate a small ball where the mass is concentrating, use the decay estimate in Corollary  \ref{space-decay} (and some additional arguments) 
to control the $L^4_t L^{12}_x$ norm in the truncated light cones emanating from that ball, and then remove those truncated
light cones to lower the energy and then apply the induction on energy hypothesis.

We turn to the details.  We fix $E \gtrsim 1$, and let $u$ be a smooth compactly supported solution to 
\eqref{nlw} on  a slab $[t_-, t_+] \times \R^3$ of energy at most $E$.  Our task is to show that
$$ \| u \|_{L^4_t L^{12}_x([t_-,t_+] \times \R^3)} \lesssim E^{O(C_2 E^{35})} (1 + M(E - c_0 E^{-33/2})).$$
We can of course assume the right-hand side is finite, thus by the standard local well-posedness theory (based largely on the Strichartz
estimate \eqref{strichartz}) any classical initial data of energy less than $E - c_0 E^{-33/2}$ will lead to
a global solution to \eqref{nlw} with $L^4_t L^{12}_x(\R \times \R^3)$ norm bounded by $M(E - c_0 E^{-33/2})$.  We shall refer to this statement
as the \emph{induction hypothesis}.

Clearly we may assume that
$$ \| u \|_{L^4_t L^{12}_x([t_-,t_+] \times \R^3)} \geq 4E^{1/12}.$$
This allows us to partition $[t_-,t_+]$ into consecutive intervals $I_1,\ldots,I_{J}$ such that
\begin{equation}\label{E12}
 E^{1/12} \leq \| u \|_{L^4_t L^{12}_x(I_j \times \R^3)} \leq 2 E^{1/12}
 \end{equation}
for all $1 \leq i \leq J$.  It then suffices to show that
$$ J \lesssim E^{O(C_2 E^{35}) } (1 + M(E - c_0 E^{-33/2}))^4.$$

Let $u_{t_+}, u_{t_-}$ be the free development of $u$ from times $t_+$ and $t_-$ respectively; one can think of these free solutions as ``scattering states'' for $u$.
Following a similar idea in \cite{tao}, we call an interval $I_j$ \emph{exceptional} if
\begin{equation}\label{exceptional-def}
 \| u_{t_+} \|_{L^4_t L^{12}_x(I_j \times \R^3)} + \| u_{t_-} \|_{L^4_t L^{12}_x(I_j \times \R^3)} 
 \geq E^{-C_2 E^{35}}
 \end{equation}
and \emph{unexceptional} otherwise.  From the Strichartz inequality \eqref{strichartz} and \eqref{energy-bound-kinetic} 
we see that there are at most $O( E^{O(C_2 E^{35}) } )$ 
exceptional intervals.  Thus it will suffice to show that any sequence $I_{j_0}, \ldots, I_{j_1}$ of consecutive
unexceptional intervals consists of at most $E^{O(C_2 E^{35}) } (1 + M(E - c_0 E^{-33/2}))^4$ intervals.  

Let us now fix $I_{j_0}, \ldots, I_{j_1}$ to be a nonempty sequence of consecutive unexceptional intervals\footnote{One may wish to think of the hypothesis that all the intervals in $I$ are unexceptional as a claim that the solution cannot ``escape'' to scatter to infinity in either time direction for the duration of the interval $I$; this can be used to establish various concentration and stability bounds on the solution here.  Compare with the irreducibility of the \emph{minimal energy blowup solutions} from \cite{gopher}, or the asymptotic smoothness of the \emph{weakly bound states} from \cite{tao:radialfocus}.  See also the concentration compactness techniques from \cite{bg}.}, and let 
$I := I_{j_0} \cup \ldots \cup I_{j_1}$.  
It will suffice to establish the bound
\begin{equation}\label{412-targ}
 \| u \|_{L^4_t L^{12}_x(I \times \R^3)} \lesssim E^{O(C_2 E^{35}) } (1 + M(E - c_0 E^{-33/2})).
 \end{equation}
From Lemma \ref{massconc} we see that for each $j_0 \leq j \leq j_1$ there exists a time $t_j \in I_j$ and a ball $B(x_j,r_j)$ 
with $r_j \lesssim E^{1/3} |I_j|$
such that
$$ \int_{B(x_j,r_j)} |u(t_j,y)|^2\ dy \geq c E^{-11/2} r_j^2$$
for some small absolute constant $c > 0$.  

Since $u$ is classical, the quantity $\frac{1}{r^2} \int_{B(x,r)} |u(t,y)|^2\ dy$ converges to $0$ as $r \to 0$ uniformly for $x \in \R^3$ and $t \in I$.  Thus we may find a minimal radius $r_0 > 0$ for which there exists a time $t_0 \in I$ and ball $B(x_0,r_0)$ with
$$  \int_{B(x_0,r_0)} |u(t_0,y)|^2\ dy = c E^{-11/2} r_0^2.$$
From the preceding discussion we see that $r_0$ exists and obeys the bound
$$ r_0 \lesssim E^{1/3} \min_{j_0 \leq j \leq j_1} |I_j|.$$

The hypotheses and conclusions so far are invariant with respect to both space translation invariance and the scaling \eqref{scaling}.  We shall now spend both of these invariances by normalising $x_0 = 0$ and $r_0 = 1$; we have two further symmetries remaining, time translation and time reversal symmetry, that we shall spend later.  We now have the useful lower bound
\begin{equation}\label{ij-lower}
|I_j| \gtrsim E^{-1/3} \hbox{ for all } j_0 \leq j \leq j_1
\end{equation}
(this basically asserts that there are no scales of interest that are significantly finer than the unit scale $1$)
and the mass bound
\begin{equation}\label{massb}
\int_{|x| \leq 1} |u(t_0,x)|^2\ dx \gtrsim E^{-11/2}.
\end{equation}
Intuitively, this means that we can remove the energy at time $t_0$ and position $B(x,1)$ from the solution and reduce the energy $E$ by a nontrivial
amount.  Unfortunately, there is a slight technical issue because the energy involves the derivative of $u$ rather than $u$ itself and so one cannot simply truncate $u$ by a rough cutoff to reduce the energy\footnote{One possible approach here would be to use equipartition of energy methods to try to force largeness of the \emph{time derivative} $|\partial_t u|$ on a ball, which would be easy to truncate.  While the author was eventually able to make this work (basically by integrating $(1+|x|^2)^{-1/2} \Box(u^2)$ in spacetime in two different ways), the bounds on the ball on which the time derivative was large ended up being exponential in $E$, which made the result comparable in nature to the somewhat simpler pigeonholing argument we give in the text.  In any event, while the exponential loss incurred here is currently dominant, there is another exponential loss arising from Corollary \ref{space-decay} and also from the induction on energy argument, and so even if the exponential loss was eliminated completely here, one would still end up with an exponential bound.}.  To resolve this we use a pigeonholing trick of Bourgain \cite{bourg.critical}.  Let $R \geq 2$ be a large radius to be chosen later.  From \eqref{hardy} and the pigeonhole principle we can find $2 < r < R$ such that
\begin{equation}\label{pigeon}
\int_{r/2 < |x| < r} |u(t_0,x)|^2\ dx \lesssim \frac{E}{\log R} r^2.
\end{equation}
Let $\chi$ be a smooth cutoff function supported on $B(0,1)$ which equals $1$ on $B(0,1/2)$, and is bounded by $1$ throughout.  
Let $\tilde u$ be the solution to \eqref{nlw} whose initial data at $t_0$ is given by
$$ \tilde u(t_0,x) := (1 - \chi(x/r)) u(t_0,x); \quad \partial_t \tilde u(t_0) := \partial_t u(t_0).$$
By comparing the energies of $u$ and $\tilde u$ at time $t_0$ we see that
\begin{align*}
E[u] - E[\tilde u] &\geq \frac{1}{6} \int_{|x| \leq 1} |u(t_0,x)|^6 - O( \frac{1}{r} \int_{r/2 < |x| < r} |u(t_0,x)| |\nabla u(t_0,x)|\ dx ) \\
&\quad - O( \int_{r/2 < |x| < r} |u(t_0,x)|^2\ dx)
\end{align*}
and hence by \eqref{massb}, \eqref{pigeon}, \eqref{energy-bound-kinetic}, and H\"older
$$ E[u] - E[\tilde u] \gtrsim E^{-33/2} - \frac{E}{\log^{1/2} R}.$$
Thus if we set $R := \exp( C E^{35} )$ for some sufficiently large constant $C$, we see that
$$ E[u] - E[\tilde u] \geq c_0 E^{-33/2}. $$
We can then apply the induction hypothesis and conclude that $\tilde u$ exists globally with the bound
$$ \| \tilde u \|_{L^4_t L^{12}_x( \R \times \R^3 )} \leq M( E - c_0 E^{-33/2} ).$$
By finite speed of propagation we conclude that
$$ \| u \|_{L^4_t L^{12}_x( \{ (t,x) \in I \times \R^3: |x| \geq R + |t-t_0| \} )} \leq M( E - c_0 E^{-33/2} )$$
and so it will suffice to show that
$$ \| u \|_{L^4_t L^{12}_x( \{ (t,x) \in I \times \R^3: |x| \leq R + |t-t_0| \} )} \lesssim E^{O(C_2 E^{35}) }.$$
By spending the time reversal symmetry (note that this can swap $u_{t_+}$ and $u_{t_-}$) it suffices to show that
$$ \| u \|_{L^4_t L^{12}_x( \{ (t,x) \in I \times \R^3: t \geq t_0; t-t_0 \leq R + |x| \} )} \lesssim E^{O(C_2 E^{35}) }.$$
We now spend the time translation symmetry to set $t_0 := R$, thus $R \in I$ and we wish to show
\begin{equation}\label{u412-targ}
\| u \|_{L^4_t L^{12}_x( \Gamma_+( I_+ ) )} \lesssim E^{O(C_2 E^{35}) }
\end{equation}
where $I_+ := I \cap [R,+\infty)$.
In doing this, we may of course assume that
$$ \| u \|_{L^4_t L^{12}_x( \Gamma_+( I_+ ) )} \geq 4 E^{1/12}.$$
From \eqref{E12}, this means that $I_+$ contains at least one full interval $I_j$ from $I$, which by \eqref{ij-lower} means that we have
the lower bound $|I_+| \gtrsim E^{-1/3}$.  (Indeed, the worst case will arise when $I_+$ is much larger than this.)

Now, from \eqref{massb} and Lemma \ref{massstable} we have
$$ \int_{|x| \leq 1} |u(t,x)|^2\ dx \gtrsim E^{-11/2}$$
provided that $|t-R| \leq c_0 E^{-13/4}$.  From H\"older's inequality we conclude that
$$ \| u \|_{L^4_t L^{12}_x( \Gamma_+([R, R + c_0 E^{-13/4}]) )} \gtrsim c_0^{1/4} E^{-57/16}.$$
Thus if we set 
\begin{equation}\label{etadef}
\eta := c_1 E^{-57/16},
\end{equation}
we have
$$ \| u \|_{L^4_t L^{12}_x( \Gamma_+([R, R + c_0 E^{-13/4}]) )} \geq 4 \eta.$$
Note that $[R, R + c_0 E^{-13/4}]$ will be contained in $I_+$ from our lower bound on $|I_+|$.

We now subdivide $I_+$ into consecutive intervals $\tilde I_1, \ldots, \tilde I_{\tilde J}$ on which
\begin{equation}\label{tilde-il}
\eta \leq  \| u \|_{L^4_t L^{12}_x( \Gamma_+(\tilde I_i) )} \leq 2 \eta 
\end{equation}
for all $1 \leq i \leq \tilde J$.  In particular, we see that 
\begin{equation}\label{first}
\tilde I_1 \subseteq [R, R + c_0 E^{-13/4}].
\end{equation}
To control the subsequent intervals, we shall use two tools.  The first is Corollary \ref{space-decay},
which will allow us to control the number of these intervals provided that none of them are too ``lacunary''.  The second tool, which we now present, gives us this desired non-lacunarity.

\begin{lemma}[Non-lacunarity estimate]\label{nonlacun}   Let $1 \leq i < \tilde J$.  Then we either have
\begin{equation}\label{tii}
|\tilde I_{i+1}| \lesssim E^{8/3} \eta^{-4} |\tilde I_i|.
\end{equation}
or
\begin{equation}\label{tii-2}
 |\tilde I_i| \geq C_2 E^{C_2 E^{35}}
\end{equation}
\end{lemma}

\begin{remark} The latter hypothesis will soon be eliminated (essentially by a ``continuity'' argument).  Intuitively, the idea is that
if $\tilde I_{i+1}$ is very much larger than $\tilde I_i$, then the portion of the solution which is concentrating in $\tilde I_i$ can ``escape''
to infinity after dispersing through $\tilde I_{i+1}$, and thus join the scattering solution $u_{t_+}$.  However, as we have eliminated all the exceptional intervals, we know that this cannot occur (unless $\tilde I_i$ is so large that it contains a very large number of unexceptional intervals).
\end{remark}

\begin{proof}  Write $\tilde I_i := [t_{i-1}, t_i]$ and $\tilde I_{i+1} = [t_i, t_{i+1}]$.  From \eqref{tilde-il}, \eqref{energy-bound-potential} and H\"older we have
$$ \| u^5 \|_{L^1_t L^2_x(\Gamma_+(\tilde I_i \cup \tilde I_{i+1}))} \lesssim E^{1/6} \eta^4$$
and hence by Duhamel \eqref{duhamel} and Strichartz \eqref{strichartz}, and finite speed of propagation, we have
$$ \| u - u_{t_{i+1}} \|_{L^4_t L^{12}_x(\Gamma_+(\tilde I_i))} \lesssim E^{1/6} \eta^4.$$
From \eqref{tilde-il} again (and the choice of $\eta$) we conclude that
$$ \| u_{t_{i+1}} \|_{L^4_t L^{12}_x(\Gamma_+(\tilde I_i))} \sim \eta.$$
This means that at least one of the statements
\begin{equation}\label{dichotomy-1}
 \| u_{t_{i+1}} - u_{t_+} \|_{L^4_t L^{12}_x(\Gamma_+(\tilde I_i))} \gtrsim \eta
\end{equation}
or
\begin{equation}\label{dichotomy-2}
 \| u_{t_+} \|_{L^4_t L^{12}_x(\Gamma_+(\tilde I_i))} \gtrsim \eta
\end{equation}
is true.

Suppose first that \eqref{dichotomy-1} holds.  From \eqref{energy-bound-kinetic} and \eqref{strichartz} we already know that
$$  \| u_{t_{i+1}} - u_{t_+} \|_{L^\infty_t L^6_x(\Gamma_+(\tilde I_i))} \lesssim E^{1/2}$$
while from the asymptotic stability estimate in Corollary \ref{ass} we have
$$ \| u_{t_{i+1}} - u_{t_+} \|_{L^\infty_t L^\infty_x(\Gamma_+(\tilde I_i))} \lesssim E^{5/6} |\tilde I_{i+1}|^{-1/2}.$$
From H\"older in space and time we conclude that
$$ \| u_{t_{i+1}} - u_{t_+} \|_{L^4_t L^{12}_x(\Gamma_+(\tilde I_i))} \lesssim E^{2/3} |\tilde I_i|^{1/4} |\tilde I_{i+1}|^{-1/4}$$
whence we obtain \eqref{tii}.

Now supppose instead that \eqref{dichotomy-2} holds.  Since none of the intervals $I_j$ that comprise $I$ are exceptional, 
we see (from the \emph{failure} of \eqref{exceptional-def}) that $\tilde I_i$ must consist of at least $\gtrsim E^{4C_2 E^{35}} / \eta^4$
intervals of the form $I_j$.
The claim \eqref{tii-2} then follows from \eqref{ij-lower}.
\end{proof}

We can now combine Lemma \ref{nonlacun} with Corollary \ref{space-decay} to obtain a key \emph{upper} bound on $I_+$.

\begin{proposition}  We have $|I_+| \leq C_2 E^{C_2 E^{35}}$.
\end{proposition}

\begin{proof} Suppose for contradiction that this estimate failed.  Let $I_{i_1}$ be the first interval for which
$|\tilde I_1 \cup \ldots \cup \tilde I_{i_1}|$ exceeds $C_2 E^{C_2 E^{35}}$.  Then by Lemma \ref{nonlacun} we have \eqref{tii} for all
$1 \leq i < i_1$.  In particular, if we let $[T_1,T_2] := \tilde I_2 \cup \ldots \cup \tilde I_{i_1-1}$, we have
$$ |T_2 - T_1| \geq |\tilde I_{i_1-1}| \gtrsim E^{-8/3} \eta^{4} |\tilde I_{i+1}|.$$
On the other hand, we have
$$ |\tilde I_1| + (T_2-T_1) + |\tilde I_{i+1}| = |\tilde I_1 \cup \ldots \cup \tilde I_{i_1}| \geq C_2 E^{C_2 E^{35}}$$
so from \eqref{first} we conclude that
$$ T_2-T_1 \gtrsim E^{C_2 E^{35} / 2 }.$$
Also from \eqref{first} and the definition of $R$ we have
$$ T_1 \lesssim R \lesssim \exp( O( E^{35} ) )$$
and thus
$$ T_2 / T_1 \gtrsim E^{C_2 E^{35}/4 }.$$
We may thus apply Corollary \ref{space-decay} (recalling the definition of $\eta$) 
and locate a time interval $[t', C_1 E^{16/3} \eta^{-8} t']$ inside $[T_1,T_2]$ such that
$$ \| u \|_{L^4_t L^{12}_x(\Gamma_+([t',C_1 E^{16/3} \eta^{-8} t']))} \leq \eta/4.$$
By \eqref{tilde-il}, this means that $[t', C_1 E^{16/3} \eta^{-8} t']$ is covered by at most two of the intervals $\tilde I_i$; from the various
inclusions of the time intervals we see that $2 \leq i < i_1$.  By the pigeonhole principle (and recalling that all of the $\tilde I_i$ lie to the right of $R$, and hence lie in the positive time axis $\R^+$) we can thus find $2 \leq i < i_1-1$ such that
 $|\tilde I_i| \lesssim C_1^{1/2} E^{8/3} |\tilde I_{i-1}|$.  But this contradicts Lemma \ref{nonlacun}.
 \end{proof}
 
From this proposition and \eqref{ij-lower} we see that $\tilde J = O(E^{O(C_2 E^{35}) })$; applying \eqref{tilde-il} we obtain \eqref{u412-targ}
as desired.  This concludes the proof of Proposition \ref{mainprop} and hence Theorem \ref{main}.

\end{document}